\documentstyle{amsart}

\newtheorem{Theorem}{Theorem}[section]

\newtheorem{Lemma}[Theorem]{Lemma}

\newtheorem{Corollary}[Theorem]{Corollary}

\newtheorem{Proposition}[Theorem]{Proposition}

\newtheorem{Remark}[Theorem]{Remark}

\newtheorem{Example}[Theorem]{Example}

\newtheorem{Definition}[Theorem]{Definition}

\begin{document}

\title{Completions of valuation rings}

\author{Steven Dale Cutkosky and Laura Ghezzi}
\thanks{Research of the first author was  partially supported by NSF}

\address{Department of Mathematics, University of Missouri,
Columbia MO 65211.} \email{cutkoskys@@math.missouri.edu \\
ghezzi@@math.missouri.edu}

\begin{abstract}
Let $k$ be a field of characteristic zero, $K$ an algebraic
function field over $k$, and $V$ a $k$-valuation ring of $K$.
Zariski's theorem of local uniformization shows that there exist
algebraic regular local rings $R_i$ with quotient field $K$ which
are dominated by $V$, and such that the direct limit $\cup R_i=V.$

We investigate the ring $T=\cup \hat R_i$. The ring $T$ is
Henselian and thus can be considered to be a ``completion'' of the
valuation ring $V$. We give an example showing that $T$ is in
general not a valuation ring. Making use of a result of Heinzer
and Sally, we give necessary and sufficient conditions for $T$ to
be a valuation ring.

The essential obstruction to $T$ being a valuation ring is the
problem of the rank of the valuation increasing upon extending the
valuation dominating a particular $R$ to a valuation dominating
its completion. In the case of rank 1 valuations, we show that
this problem can be handled in a very satisfactory way.

Finally, suppose that $K^*$ is a finite algebraic extension of $K$
and $V^*$ is a rank 1 $k$-valuation ring of $K^*$ such that
$V=V^*\cap K$. We obtain a relative local uniformization theorem
for the extension $K^*$ of $K$, that generalizes previous results
of Cutkosky and Piltant.

\end{abstract}

\maketitle

\section{Introduction}
Suppose that $k$ is a field of characteristic zero, $K$ is an
algebraic function field over $k$, and $V$ is a $k$-valuation ring
of $K$ (that is, $k\subset V$ and the quotient field of $V$ is
$K$).  Zariski's theorem of local uniformization \cite{Z} shows
that there exist algebraic regular local rings $R_i$ with quotient
field $K$ which are dominated by $V$, and such that the direct
limit
$$\cup R_i=V.$$
Now suppose that $K^*$ is a finite algebraic extension of $K$ and
$V^*$ is a $k$-valuation ring of $K^*$ such that $V=V^*\cap K$.
Let $\Gamma^*$ be the value group of $V^*$, and $\Gamma$ be the
value group of $V$.

The first author has shown with Olivier Piltant in \cite{Dale} that a relative local uniformization theorem holds for the  extension $K^*$ of $K$, which gives the strongest
possible generalization of the classical ramification theory of Dedekind domains to
general valuations. The following theorem is a summary of some of the conclusions of
Theorem 6.3 \cite{Dale}.

\begin{Theorem}\label{TheoremI1} Suppose that assumptions and notations are as above.
Let $k'$ be an algebraic closure of $V^*/m_{V^*}$. Then there
exist a directed system of algebraic regular local rings $S_i$
with quotient field $K^*$ which are dominated by $V^*$, and a
directed system of algebraic normal local rings $R_i$ with
quotient field $K$ which are dominated  by $V$ such that
\begin{enumerate}
\item $\cup S_i=  V^*$ and $\cup R_i=V$.

 \item $S_i$ is a
localization at a maximal ideal of the integral closure of $R_i$
 in $K^*$ for all $i$.
\item There exist  actions of $\Gamma^*/\Gamma$ on $\hat
S_i\otimes_{S_i/m_{S_i}}k'$ which are compatible with the directed
system so that
$$
(\hat S_i\otimes_{S_i/m_{S_i}}k')^{\Gamma^*/\Gamma}\cong \hat
R_i\otimes_{R_i/m_{R_i}}k'.
$$
\end{enumerate}
\end{Theorem}
It was shown by an example of Abhyankar \cite{Ab4}  that  it is in
general  not possible to find an algebraic regular local ring $S$
with quotient field $K^*$ which is dominated by $V^*$ such that
there exists an algebraic regular local ring $R$ with quotient
field $K$ such that $S$ is a localization of the integral closure
of $R$ in $K^*$.   The fact (proven in \cite{Dale} and
\cite{Dale4}) that normal local rings $R$ exist satisfying this
property proves the ``local weak simultaneous resolution conjecture'' of
Abhyankar, posed by Abhyankar in \cite{Ab4} and \cite{Ab5}. The
$R_i$ found in the proof of Theorem \ref{TheoremI1} in fact have
toric singularities. This is reflected in the fact stated above
that their completions are abelian quotient singularities.

From this theorem we obtain the following.
\begin{Theorem}\label{TheoremI2} Let notations be as in Theorem  \ref{TheoremI1}.
Let $U^*=\cup \hat S_i\otimes_{S_i/m_{S_i}}k'$ and let $U=\cup
\hat R_i\otimes_{R_i/m_{R_i}}k'$. Then $U^*$ and $U$ are Henselian
normal domains, and   $Q(U^*)$ is a finite  Galois extension of
$Q(U)$ with Galois group $\Gamma^*/\Gamma$.
\end{Theorem}

We give a proof of Theorem \ref{TheoremI2} in Section \ref{Section1}.

In this paper we compare the ``completion'' of Theorem \ref{TheoremI2} with other
notions of completion of a valuation ring (\cite{K}, \cite{Kr2}, \cite{KPR}, \cite{HS}, \cite{S}, \cite{T}).

Let us briefly allow $k$ to be an arbitrary field. We summarize
some of the results of Section \ref{Section2}. Suppose that
$\{R_i\}$ is a directed system of normal algebraic local rings
which are dominated by $V$, and such that $\cup R_i=V$.  The ring
$T=\cup \hat R_i$ does not depend on our choice of $\{R_i\}$ whose
union is $V$ (Lemma \ref{LemmaA5}), and is Henselian (Proposition
\ref{PropV1}).
 Thus $T$ can be considered
to be a  ``completion'' of the valuation ring $V$. We give an
example showing that $T$ is in general not a valuation ring, and
we show that  $T$ is itself a valuation ring if and only if for
each $i$ there exists a unique valuation ring $V_i$ with quotient
field $K_i$ (where $K_i$ is the quotient field of $\hat R_i$)
which dominates $V$ and $\hat R_i$. (Theorem \ref{Theorem1}). We
make use of a theorem of Heinzer and Sally \cite{HS} on the
uniqueness of extensions of valuations dominating a local ring to
their completion in proving this result.

We give an example (Example \ref{ExampleFE25}) showing that even if $V$ and $T$ are rank 1 valuation rings, then $T$ is
in general not complete and in particular is not a maximal immediate extension, as defined in
\cite{Kr2} and \cite{K}.

The essential obstruction to $T$ being a valuation ring is the
problem of the rank of the valuation increasing upon extending the
valuation dominating a particular $R$ to a valuation dominating
its completion (Corollary to Theorem \ref{Theorem1}). In the case
of rank 1 valuations, this problem can be handled in a very
satisfactory way, and (in characteristic zero) we will obtain a
good valuation theoretic explanation of Theorem \ref{TheoremI2}.

In Section \ref{Section3}, we define the prime ideal $p(\hat
R)_{\infty}$ of elements of infinite value of the completion $\hat
R$ of an algebraic local ring $R$ dominated by a rank 1 valuation $V$.
This prime has previously been defined and considered in
\cite{Dale2} and \cite{T}, as well as by Spivakovsky. The
essential point here is that there is a unique extension of the
valuation ring $V$ to a valuation ring of  the quotient field of
$\hat R/p(\hat R)_{\infty}$ which dominates $\hat R/p(\hat
R)_{\infty}$.  We  conclude that there is a unique valuation ring
$\overline V$ of the quotient field of the ring $\overline T=\cup
\hat R_i/p(\hat R_i)_{\infty}$ which contains $T$.

In the case when $V$ has rank greater than 1 there is no natural
ideal in $\hat R$ which contains the obstruction to the jumping of
the rank of an extension of $V$ to $\hat R$, although this
obstruction is obtained in  a series of prime ideals in quotient
rings of $\hat R$.

For the remainder of this introduction we assume that $k$ has characteristic zero, and
$V$ has rank 1.
We prove that $\overline T$ is in fact a valuation ring  in Theorem \ref{TheoremA6}, and that
$(\overline T, Q(\overline T))$ is an immediate Henselian extension of  $(V,K)$ in Theorem
\ref{TheoremA4}.

We further show (in Theorem \ref{TheoremA4}) that we can choose
our system of regular local rings $R_i$ so that each $\hat
R_i/p(\hat R_i)_{\infty}$ is a regular local ring. The main new
technical result used in this statement is Theorem
\ref{TheoremFE6}, which shows that we can simultaneously resolve
the primes of infinite value in a finite extension. In this case
the ``finite extension'' is just the identity, but we will need
this more general result later.

We now turn to an analysis of  our finite extension $K^*$ over
$K$, in the case when $V^*$ (and $V=V^*\cap K$) are rank 1
valuation rings and $k$ has characteristic zero.
  We make essential use of  Theorem \ref{TheoremFE6} (on simultaneous resolution of the primes of
infinite value).
We obtain in Theorems \ref{Theorem8} and \ref{TheoremFE15} a generalization of
Theorem 5.1 \cite{Dale2} and Theorem \ref{TheoremI1} (Theorem 6.3 \cite{Dale}) in this context.
   Let $k'$ be an algebraic closure of $V^*/m_{V^*}$. We find a system of regular
   local rings $S_i$ whose union is $V^*$,
and a system of normal local rings $R_i$ whose union is $V$, such
that for all $i$, $S_i$ is a localization at a maximal ideal of
the integral closure of $R_i$ in $K^*$. If $q_i$ and $p_i$ are the
respective primes of infinite value, then $\hat S_i/q_i$ is a
regular local ring and $\hat R_i/p_i$ is a normal local ring with
toric singularities. There are compatible actions of
$\Gamma^*/\Gamma$ on ($\hat S_i/q_i)\otimes_{S_i/m_{S_i}}k'$ such
that
$$
((\hat S_i/q_i)\otimes_{S_i/m_{S_i}}k')^{\Gamma^*/\Gamma}\cong
(\hat R_i/p_i) \otimes_{R_i/m_{R_i}}k'.
$$

From  Theorems \ref{Theorem8} and \ref{TheoremFE15} we obtain the following.
\begin{Theorem}\label{TheoremI3} Let notations be as in Theorem  \ref{TheoremI1}.
Assume that $V^*$ (and  $V=V^*\cap K$) have rank 1, and that
$k=V^*/m_{V^*}$ is algebraically closed of characteristic zero.
Then there exist directed systems of algebraic local rings $R_i$ and
$S_i$ satisfying the conclusions of Theorem \ref{TheoremI1} and such that
$\hat R_i/p_i$ and $\hat S_i/q_i$,
where $q_i$ and $p_i$ are the primes of
elements of infinite value, satisfy the conclusions of the above paragraph. Let $\overline U^*=\cup \hat S_i/q_i$ and let $\overline
U=\cup\hat R_i/p_i$. Then $\overline U^*$ and $\overline U$
are Henselian valuation rings, such that $(\overline U^*,
Q(\overline U^*))$ and $(\overline U, Q(\overline U))$ are
immediate extensions of $(V^*,K^*)$ and $(V,K)$ respectively, and
$Q(\overline U^*)$ is a finite  Galois extension of $Q(\overline
U)$ with Galois group $\Gamma^*/\Gamma$.
\end{Theorem}

We give the proof of Theorem \ref{TheoremI3} in Section \ref{Section6}.

\section{Notations}
We will denote the
maximal ideal of a local ring $R$ by $m_R$ or $m(R)$.
We will denote the  quotient field of a domain $R$ by $Q(R)$.
Suppose that $R\subset S$ is an inclusion of local rings. We will say that $R$ dominates
$S$ if $m_S\cap R=m_R$.
Suppose that $K$ is an algebraic function field over a field $k$. We will say that
a subring  $R$ of $K$ is algebraic if $R$ is essentially of finite type over $k$.
Suppose that $K^*$ is a finite extension of an algebraic function field $K$,
 $R$ is a local ring with quotient field $K$ and $S$ is a local ring with quotient field $K^*$. We will say
that $S$ lies over $R$ and $R$ lies below $S$ if $S$ is a localization at a maximal ideal of the integral closure of $R$ in $K^*$.
If $R$ is a local ring, $\hat R$ will denote the completion of $R$ at its maximal ideal.

Good introductions to the valuation theory which we require in this paper can be found  in Chapter VI of \cite{ZS} and in \cite{Ab}.
A valuation $\nu$ of $K$ will be called a $k$-valuation if $\nu(k)=0$.
We will denote by $V_{\nu}$ the associated valuation ring, which necessarily contains $k$.
A valuation ring
$V$ of $K$ will be called a $k$-valuation ring if $k\subset V$.
The residue field
$V/m_V$ of a valuation ring $V$ will be denoted by $k(\nu)$. The value group of a valuation $\nu$ with valuation ring $V$ will be denoted by $\Gamma_{\nu}$ or $\Gamma_V$.
If $R$ is a subring of $V_{\nu}$ then the center of $\nu$ (the center of $V_{\nu}$)
on $R$ is the prime ideal $R\cap m_{V_{\nu}}$.

Suppose that $R$ is a local domain. A monoidal transform $R\rightarrow R_1$ is a
birational extension of local domains such that $R_1=R[\frac{P}{x}]_m$ where $P$ is
a regular prime ideal of $R$, $0\ne x\in P$ and $m$ is a prime ideal of $R[\frac{P}{x}]$
such that $m\cap R=m_R$. $R\rightarrow R_1$ is called a quadratic transform  if $P=m_R$.

If $R$ is regular, and $R\rightarrow R_1$ is a monoidal transform,
then there exists a regular system of parameters $(x_1,\ldots, x_n)$ in
$R$ and $r\le n$ such that
$$
R_1=R\left[\frac{x_2}{x_1},\ldots,\frac{x_r}{x_1}\right]_m.
$$

Suppose that $\nu$ is a valuation of the quotient field $R$ with valuation ring $V_{\nu}$
which dominates $R$. Then $R\rightarrow R_1$ is a monoidal transform along $\nu$
(along $V_{\nu}$) if $\nu$ dominates $R_1$.

\section{Completion of relative local uniformization}\label{Section1}
We now give the proof of Theorem  \ref{TheoremI2}.

Let $\overline S_i=\hat S_i\otimes_{S_i/m_{S_i}}k'$, $\overline
R_i=\hat R_i\otimes_{R_i/m_{R_i}}k'$ for $i\in I$ and let
$G=\Gamma^*/\Gamma$. $U^*$ is the directed union of the $\overline
S_i$ and $U$ is the directed union of the $\overline R_i$. The
fact that $R_i$ is a normal domain implies  $\hat R_i$ is a normal
domain (Chapter VIII,  Section 13, Theorem 32 \cite{ZS}) and thus
$\overline R_i$ is a normal local domain (Proposition IV.6.7.4
\cite{EGA}). Thus $U$ (and $U^*$) are normal domains. The fact
that $U$ and $U^*$ are Henselian follows from the proof of
Proposition \ref{PropV1}. The action of $G$ on $U^*$ extends to an
action on $Q(U^*)$. Suppose that $h\in U^*$ and $\sigma(h)=h$ for
all $\sigma\in G$. There exists $i$ such that $h\in Q(\overline
S_i)$. Since $\overline S_i$ is finite over $\overline R_i$ and
$\overline S_i^{G}=\overline R_i$, it follows that $Q(\overline
S_i)^G=Q(\overline R_i)$.  Thus $h\in Q(U)$. We conclude that
$Q(U^*)^G=Q(U)$, so that $Q(U^*)$ is a finite Galois extension of
$Q(U)$ (c.f. Theorem V.2.15 \cite{Hu}).

\begin{Remark} The statement that $Q(U^*)$ is finite over $Q(U)$ can be seen directly from the
fact that the minimal polynomial of each $Q(\hat S_i)$ over
$Q(\hat R_i)$ is a factor of the minimal polynomial of an
appropriate primitive element of $K^*$ over $K$ $($by Proposition
$1$ \cite{Ab3}$)$.
\end{Remark}

The fact that the extension considered in  Theorem \ref{TheoremI2} is Galois, even
when the original field extension $K^*/K$ is not, is a condition that can be easily
seen in the case when $U$ and $U^*$ are valuation rings, as the first author realized
 with Franz-Viktor Kuhlmann in a discussion.

\begin{Theorem}\label{TheoremFE17} Suppose that  $V$ is a Henselian valuation ring
of a field $K$, such that $V$ contains an algebraically closed field $k$
 of characteristic zero, with $k\cong V/m_V$.
If $L$ is a finite extension of $K$, then there is a unique
valuation ring $W$ of $L$ such that $W$ dominates $V$, and $L$ is
Galois over $K$ with Galois group $\Gamma_W/\Gamma_V$, where
$\Gamma_W$ and $\Gamma_V$ are the respective value groups.
\end{Theorem}
\begin{pf}
Let $J$ be a finite Galois extension of $K$ which contains $L$.
Let $G$ be the Galois group of $J$ over $K$. Since $V$ is
Henselian, there exists a unique valuation ring $U$ of $J$ such
that $U$ dominates $V$ ((16.4), (16.6) \cite{E}). Thus the
splitting group $G^s(U/V)=G$ by Proposition 1.46 \cite{Ab}. We
have $U/m_U=V/m_V=k$ since $k$ is algebraically closed. Thus the
inertia group $G^i(U/V)=G^s(U/V)=G$ by Theorem 1.48 \cite{Ab}.
Finally,
$$
G=G^i(U/V)\cong \Gamma_U/\Gamma_V
$$
by Theorem 3 \cite{Kr} or Chapter VI, Section 12, Corollary \cite{ZS}.

Since $G$ is abelian, all intermediate subfields of $J$ are Galois
over $K$. Thus $L$ is Galois over $K$, and the Galois group of $L$
over $K$ is $\Gamma_W/\Gamma_V$.
\end{pf}

\section{Completions of valuation rings}\label{Section2}

Suppose that $K$ is an algebraic function field over a field $k$, and $V$ is a valuation
ring of $K$ with maximal ideal $m_V$ and value group $\Gamma$.  Suppose that
 $\{R_i\mid i\in I\}$ is a directed system of normal local rings such that
\begin{enumerate}
\item[(a)]  $V=\cup_{i\in I}R_i$. \item[(b)] $I$ has a minimum
$0$. \item[(c)] Each $R_i$ is essentially of finite type over $k$
and has quotient field $K$. \item[(d)] If $i<j$ then $R_j$
dominates $R_i$.
\end{enumerate}
 Let $K_i$ be
the quotient field of $\hat R_i$.  By Zariski's subspace theorem
((10.13) \cite{Ab2}) we have natural inclusions $\hat
R_i\rightarrow \hat R_j$ if $i<j$, and $\{\hat R_i\mid i\in I\}$
is a directed system of normal local rings (Scholie 7.8.3
\cite{EGA}). Let
  $$
T=\cup_{i\in I} \hat R_i
$$
 and
$K_{\infty}=\cup K_i$. $T$ is a normal domain with quotient field
$K_{\infty}$ and maximal ideal $m_T=\cup m_{\hat R_i}$.

Given a valuation ring $V$ as above, there exists a directed system of normal local rings
$\{R_i\}$ whose union is $V$.  A particular construction is as follows.
We take $R_0$ to be any normal local ring which is dominated by $V$.
If $m\in{\bf N}$ and $f_1,\ldots,f_m\in V$ we set $i=(f_1,\ldots,f_m)$ and let
$R_i$ be the localization of the normalization of $R_0[f_1,\ldots,f_m]$ which is dominated by
$V$.

\begin{Lemma}\label{LemmaA5}
The ring $T=\cup_{i\in I}\hat{R_i}$ is independent of choice of
directed system $\{R_i\mid i\in I\}$ satisfying $(a)$,$(b)$,$(c)$,
and $(d)$.
\end{Lemma}
\begin{pf} Let $J$ be a partially
ordered set, and let $\{S_j|\  j\in J\}$ be a collection of
algebraic local rings with quotient field $K$, such that $\{S_j\}$
satisfies (a),(b),(c) and (d). We show that $\cup_{i\in I}\hat
R_i=\cup_{j\in J}\hat {S_j}$.

Let $i\in I$. Since $R_i$ is essentially of finite type over $k$
and dominated by $V$, there exist $f_1,\dots,f_m \in V$ such that
$R_i=k[f_1,\dots,f_m]_{ m_V\cap k[f_1,\dots,f_m]}$.
Since $V=\cup_{j\in J} S_j$ and $J$ is directed, there exists
$j\in J$ such that $f_1,\dots,f_m\in S_j$, and so
$k[f_1,\dots,f_m]_{ m_V\cap k[f_1,\dots,f_m]}\subset
S_j$, since $S_j$ is dominated by $V$. Hence $R_i\subset S_j$.
There is then a natural inclusion $\hat R_i\subset \hat S_j$, and thus
$\cup_{i\in I}\hat{R}_i\subset\cup_{j\in
J}\hat {S}_j$. The other inclusion is proven in the same
way.
\end{pf}

\begin{Proposition}\label{PropV1}
The ring $T=\cup_{i\in I} \hat R_i$ is Henselian.
\end{Proposition}

\begin{pf} Let $F\in T[x]$ and $\phi_1, \phi_2 \in T/m_T[x]$ be
monic polynomials such that $\phi_1$ and $\phi_2$ are relatively
prime and $\bar F=\pi(F)=\phi_1 \phi_2$, where $\pi:T\rightarrow
T/m_T$ is the natural projection. We need to show that there exist
monic polynomials $F_1, F_2 \in T[x]$ such that $\pi(F_1)=\phi_1$,
$\pi(F_2)=\phi_2$, and $F=F_1 F_2$.

Since both $T=\cup \hat R_i$ and $T/m_T=\cup \hat R_i/m_{\hat R_i}$ are directed unions,
there exists $c\in I$ such that $F\in \hat R_c[x]$ and $\phi_1,\phi_2\in\hat R_c/m_{\hat R_c}$.

 Since $\hat R_c$ is complete, there exist monic
polynomials $F_1, F_2 \in \hat R_c[x]$ such that
$\pi'(F_1)=\phi_1$, $\pi'(F_2)=\phi_2$, and $F=F_1 F_2$. Since
$F_1, F_2 \in T[x]$, $F_1$ and $F_2$ are the desired polynomials.
\end{pf}

\begin{Example} In general, $T$ is not a valuation ring.
\end{Example}

\begin{pf}
Let $k$ be a field and $R=k[x,y,z]_{(x,y,z)}$. We will define a valuation $\nu$ on
$K=Q(R)$ which dominates $R$.  Let
$p(t)=\sum_{i=1}^{\infty}a_it^i\in k[[t]]$ be a transcendental series. If
$f\in k[x,y,z]$ write $f=z^rg(x,y,z)$ where $g(x,y,0)\ne 0$. We
define
$$
\nu(f)=(r,\text{ ord }g(t,p(t),0)).
$$
The value group of $\nu$ is ${\bf Z}\oplus{\bf Z}$ with the
lexicographic order. Let $V$ be the valuation ring of $\nu$ and
let $R_i$ be a directed system of regular local rings satisfying
(a), (b), (c) and (d) with $R_0=R$. Such a system exists by
Theorem \ref{TheoremI1}. In particular, $V=\cup R_i$.

We will suppose that $T$ is a valuation ring and derive a
contradiction. Let $\overline\nu$ be an extension of $\nu$ to $Q(T)$
 such
that $T$ is the valuation ring of $\overline \nu$.  There is a natural
embedding of value groups
$$
\Gamma_{\nu}={\bf Z}\oplus{\bf Z}\subset \Gamma_{\overline \nu}.
$$
$\frac{z}{y-p(x)}\in Q(T)$. If $\overline\nu(\frac{z}{y-p(x)})\ge 0$
let $f=\frac{z}{y-p(x)}$. If $\overline\nu(\frac{z}{y-p(x)})<0$
let $f=\frac{y-p(x)}{z}$. By construction $\overline\nu(f)\ge 0$
so that $f\in T$ (by our assumption that $T$ is a valuation ring). Thus there exists $\hat R_i$ such that $f\in\hat
R_i$. By our hypothesis, $R_i$ dominates $R$.  Since $R_i$ is essentially of finite type over
$R$, $Q(R)=Q(R_i)$ and $R$ is a UFD, there exists an ideal $I\subset R$  of height
$\ge 2$ and $a\in I$ with
$$
\nu(a)=\text{min}\{\nu(h)\mid h\in I\}
$$
such that  $R[\frac{I}{a}]\subset V$ and
$$
R_i=R[\frac{I}{a}]_{m_V\cap R[\frac{I}{a}]}
$$
Since $I\not\subset (z)$, we have that $\nu(a)\in\{0\}\oplus{\bf
Z}$. Let $R'=\hat R[\frac{I}{a}]_{m_T\cap \hat R[\frac{I}{a}]}$.
$$
f\in Q(\hat R)\cap \hat R_i=Q(R')\cap \hat R'=R'.
$$
For the last equality, c.f. Lemma 2.1 \cite{Dale2}.

We can thus write $f=\frac{g}{h}$ with $g,h\in\hat R[\frac{I}{a}]$
and $h\not\in m_T\cap\hat R[\frac{I}{a}]$. Thus
$\overline\nu(h)=0$. There exists $n\in{\bf N}$ such that
$a^ng=g_0$, $a^nh=h_0$ and $g_0,h_0\in \hat R$. Thus
$\overline\nu(h_0)=n\nu(a)+\overline\nu(h)=(0,m)$ for some
$m\in{\bf N}$.

If $f=\frac{z}{y-p(x)}$, we have $zh_0=(y-p(x))g_0$ implies
$y-p(x)$ divides $h_0$ in $\hat R$. But
$$
y-p(x)=
(y-\sum_{i=1}^{m+1}a_ix^i)-x^{m+2}(\sum_{i=m+2}^{\infty}a_ix^{i-m-2}).
$$
Thus $\overline\nu(y-p(x))>(0,m)$, which is a contradiction.

If $f=\frac{y-p(x)}{z}$, then $(y-p(x))h_0=zg_0$ implies $z$
divides $h_0$ in $\hat R$.

 This is a contradiction since $\nu(z)=(1,0)>(0,m)=\overline\nu(h_0)$.

\end{pf}

As an extension of the above example, we construct valuations $\nu_1$, $\nu_2$
of $Q(\hat R)$ which extend $\nu$ and dominate $\hat R$ such that $\nu_1$ is an immediate
extensions of $\nu$ ($k(\nu)=k(\nu_1)$ and $\Gamma_{\nu}=\Gamma_{\nu_1}$) but $\nu_2$
is a rank 3 valuation, of higher rank than $\nu$.

We first define $\nu_1$. For $f\in \hat R$, write
\begin{equation}\label{eqnu1}
f(x,y,z)=\sum a_{ij}(y-p(x))^iz^j
\end{equation}
with $a_{ij}\in k[[x]]$. Let
$$
a=\text{min}\{i+j\mid a_{ij}(x)\ne 0\},
$$
$$
b=\text{min}\{\text{ord}(a_{ij}(x))\mid i+j=a\}.
$$
Set $\nu_1(f)=(a,b)\in{\bf Z}\oplus{\bf Z}$.
$\nu_1$ defines a valuation on $Q(\hat R)$ which dominates $\hat R$.
We will now verify that $\nu_1$ extends $\nu$.  For $f\in k[x,y]$, write
$$
f=\sum b_j(x,y)z^j
$$
with each $b_j\in k[x,y]$. Suppose that $\nu(f)=(r,s)$. Comparing with (\ref{eqnu1}), we see that for all $j$,
$$
b_{j}(x,y)=\sum_{i=0}^{\infty}a_{ij}(x)(y-p(x))^i.
$$
Thus we either have that $b_j=0$ or $a_{0j}\ne 0$, since $y-p(x)\not\,\mid b_j(x,y)$.
Thus
$$
a=\text{min}\{j\mid b_j(x,y)\ne 0\}=r.
$$
Set $g=\frac{f}{z^r}$.  We have
$$
b= \text{ord}(a_{0r}(x))=\text{ord}(g(x,p(x)))=s.
$$
Thus $\nu_1(f)=\nu(f)$. It follows that $\nu_1$ is an immediate extension of $\nu$.

Now we define the extension $\nu_2$.  For $f\in \hat R$,
write $f=z^{\alpha}g(x,y,z)$ where
\linebreak $z\not\,\mid g(x,y,z)$. Write
$g(x,y,0)=(y-p(x))^{\beta}h(x,y)$ where $y-p(x)\not\,\mid g(x,y,0)$.
Set $\gamma=\text{ord}(h(x,p(x)))$.
Define
$\nu_2(f)=(\alpha,\beta,\gamma)\in{\bf Z}^3$, where ${\bf Z}^3$ has the lexicographic order.
$\nu_2$ extends to a valuation of $Q(\hat R)$ which dominates $\hat R$, and such that
$\nu_2$ extends $\nu$.

\begin{Theorem}\label{Theorem1} $T$ is  a valuation ring if and only if for all $i$, there exists a unique
valuation ring $V_i$ with quotient field $K_i$ which dominates $V$
and $\hat R_i$.
\end{Theorem}

\begin{pf} Suppose that for all $i$, there exists a unique
valuation ring $V_i$ with quotient field $K_i$ which dominates $V$
and $\hat R_i$. By Proposition 4.1 \cite{HS}
$$
V_i=\hat R_i[V]_{m_V\hat R_i[V]}\subset K_i
$$
is a valuation ring for all $i$. Let
$$
V_{\infty}=\cup V_i=(\cup \hat R_i)[V]_{m_V(\cup \hat
R_i)[V]}\subset K_{\infty}
$$
where $K_{\infty}$ is the quotient field of $T$. $f\in K_{\infty}$
implies $f\in K_i$ for some $i$ which implies $f\in V_i$ or
$\frac{1}{f} \in V_i$. Thus $f\in V_{\infty}$ or $\frac{1}{f}\in
V_{\infty}$. Thus $V_{\infty}$ is a valuation ring.  $V\subset
\cup \hat R_i$  implies $V_{\infty}=\cup
\hat R_i=T$, and we conclude that $T$ is a valuation ring.

Now suppose that $T$ is  a valuation ring. Let $m=\cup m_{R_i}\hat
R_i$ be the maximal ideal of $T$. Without loss of generality, we
may assume that $V/m_V$ is algebraic over $k$. For if this is not
the case, we can replace $k$ with a rational function field $k'$
over $k$ contained in all of the $R_i$ such that $V/m_V$ is
algebraic over $k'$.

Suppose that for some index $i$, $V_1$ is a valuation ring with
quotient field $K_i$ which dominates $\hat R_i$ and $V$. We will
show that there exists a valuation ring $W_1$ with quotient field
$K_{\infty}$ which dominates $T$ and such that $W_1\cap K_i=V_1$.

Consider the domain $A=T[V_1]\subset K_{\infty}$. Let $I\subset A$
be the ideal $I=(m+m_{V_1})A$.  We will first establish that
$1\not\in I$. If it were true that $1\in I$, then there would
exist an index $j>i$ such that $1\in (m_{R_j}+m_{V_1})\hat
R_j[V_1]$. Since $R_i\rightarrow R_j$ is birational, there exists
an ideal $b\subset R_i$ and $x\in b$ such that $\frac{b}{x}\subset
V$ and $R_j$ is a localization of $R_i[\frac{b}{x}]$ at a maximal
ideal. Thus $\hat R_i[\frac{b}{x}]\subset V_1$ since
$V\subset V_1$. Let $\overline \nu$ be a valuation of $K_i$ which
has $V_1$ for its valuation ring. We have $\overline\nu(\hat
R_i[\frac{b}{x}])\ge 0$ and $\overline\nu(m_{R_j})>0$ since
$m_{R_j}=m_{V_1}\cap R_j$. Thus $m_{R_j}\hat R_i[\frac{b}{x}]\subset m_{V_1}$, and $V_1$ dominates $\hat
R_i[\frac{b}{x}]_{m_{R_j}\hat R_i[\frac{b}{x}]}$. Thus there
exists a valuation ring  $U_1$ of $K_j$ such that $U_1\cap
K_i=V_1$ and $U_1$ dominates $\hat R_j$ (as follows from  page 177
of \cite{HS}). $(m_{R_j}+m_{V_1})\hat R_j[V_1]\subset m_{U_1}$
implies $1\not\in (m_{R_j}+m_{V_1})\hat R_j[V_1]$. Thus we have a
contradiction, and $1\not\in I$.

Let $a$ be a prime ideal in $A$ which contains $I$. Suppose that
$h\in A/a$. There exists an index $j$ such that we can write $h$
as a class $h=[\sum f_rg_r]$ with $f_r\in V_1$ and $g_r\in \hat
R_j$. We have natural inclusions $\hat R_j/m_{\hat R_j}\rightarrow
A/a$ and $V_1/m_{V_1}\rightarrow A/a$ such that $h$ is the image
of the induced map $\hat R_j/m_{\hat
R_j}\otimes_kV_1/m_{V_1}\rightarrow A/a$. Thus $h$ is algebraic
over $V_1/m_{V_1}$, since $\hat R_j/m_{\hat R_j}$ is finite over
$k$. We conclude that $A/a$ is algebraic over $V_1/m_{V_1}$.

There exists a valuation ring $W_1$ which contains $A$ such that
$m_{W_1}\cap A=a$ and $W_1/m_{W_1}$ is algebraic over $A/a$, by
Corollary 3 to Theorem 5' of Section 4, Chapter VI \cite{ZS}.
Let $\overline W_1=W_1\cap K_i$. $\overline W_1$ contains $V_1$
and $m_{V_1}\subset m_{\overline W_1}$. But $\overline
W_1/m_{\overline W_1}$ is algebraic over $V_1/m_{V_1}$ since
$W_1/m_{W_1}$ is algebraic over $V_1/m_{V_1}$.  Thus
$V_1=\overline W_1$ by Theorem 2 of Section 3, Chapter VI
\cite{ZS} and Corollary 1 to Theorem 5 of Section 4, Chapter VI
\cite{ZS}.

We have thus proved the existence of an extension $W_1$ of $V_1$
to $K_{\infty}$ which restricts to $V_1$ and dominates $T$.

Continuing with the proof of the theorem, suppose that for some
index $i$, the extension of $V$ to $K_i$ which dominates $\hat
R_i$ is not unique. There are then
 extensions $V_1$ and $V_2$ of $V$ to
$K_i$ which dominate $\hat R_i$  such that $V_1\not\subset V_2$
and  $V_2\not\subset V_1$.

We have shown that there then exist valuation rings $W_1$ and
$W_2$ of $K_{\infty}$ such that $T\subset W_1\cap W_2$, $W_1\cap
K_i=V_1$ and $W_2\cap K_i=V_2$. Thus $W_1\not\subset W_2$ and
$W_2\not\subset W_1$. But this is impossible since $T$ is a
valuation ring of $K_{\infty}$, by Theorem 3, Section 3, Chapter
VI \cite{ZS}.
\end{pf}

\begin{Corollary} $T$ is a valuation ring   if for all $i$ there does not
exist an extension of $V$ to $K_i$ which dominates $\hat R_i$ of
higher rank than the rank of $V$.
\end{Corollary}

\begin{pf} This follows from Theorem \ref{Theorem1} and the remark on page 181 of
\cite{HS} which shows that if the extension of $V$ to $K_i$ which
dominates $\hat R_i$ is not unique then there must be an extension
of higher rank.
\end{pf}

The converse to the above corollary is false, as is seen by the
following simple example. Let $p(t)\in k[[t]]$ be a transcendental
power series with constant term zero. Consider the rank 1 discrete valuation $\nu$ on
$k(x,y)$ defined by the embedding of $k$-algebras
$$
k(x,y)\rightarrow k\ll t \gg
$$
generated by $x=t, y=p(t)$ where $k\ll t \gg$ denotes the quotient field of $k[[t]]$.
$\nu$ dominates $R=k[x,y]_{(x,y)}$.  The valuation ring $V$ of $\nu$ extends uniquely to a
rank 2 valuation ring of the quotient field of $k[[x,y ]]$
which dominates $k[[x,y]]$. Furthermore, the construction gives a unique
extension of $V$ to a rank 2 valuation ring which dominates $Q(\hat S)$
for any algebraic normal local ring $S$ of $k(x,y)$ such that $V$
dominates $S$ and $S$ dominates $R$.

It follows from Theorem \ref{Theorem1} that the examples in \cite{HS}
of valuation rings dominating regular local rings $R$ which do not have unique extensions in $Q(\hat R)$ dominating $\hat R$ generate examples
where $T$ is not a valuation ring.

\section{The prime ideal of elements of infinite value}\label{Section3}
We will assume in this section that $V$ has rank 1, that is, the value group of $V$
is a (possibly  nondiscrete) subgroup $\Gamma$ of ${\bf R}$. Other notations and assumptions will
be as in Section \ref{Section2}.

\begin {Lemma}\label{LemmaA1} Suppose that $V$ has rank 1 and that $R$ is an algebraic normal local ring of $K$
such that $V$ dominates $R$ and $f\in\hat R$. Then one of the following must hold.
\begin{enumerate}
\item There exists $\rho\in\Gamma$ such that if $\{f_n\}$ is any Cauchy sequence in $R$
which converges to $f$, then $\nu(f_n)=\rho$ for all $n\gg 0$.
\item If $\rho\in \Gamma$ and if $\{f_n\}$ is any Cauchy sequence in $R$ which converges to $f$,
then $\nu(f_n)>\rho$ for $n\gg 0$.
\end{enumerate}
\end{Lemma}

\begin{pf} We first argue that (1) or (2) must hold for a fixed Cauchy sequence $\{f_n\}$
in $R$ which converges to $f$.  Suppose that (2) doesn't hold.
Then there exists $\overline \rho\in\Gamma$ such that given
$n_0\in{\bf N}$, there exists $\overline n>n_0$ such that
$\nu(f_{\overline n})\le \overline \rho$. Let $i\in{\bf N}$ be
such that $i\nu(m_R)>\overline \rho$, and let $n_0$ be such that
$f_m-f_n\in m_R^i$ if $m,n>n_0$. There exists $\overline n>n_0$
such that $\nu(f_{\overline n})\le\overline \rho$. Then
$\nu(f_n)=\nu(f_{\overline n})$ if $n> n_0$, so (1) holds for
$\{f_n\}$.

If $\{f_n\}$ and $\{g_l\}$ are two distinct Cauchy sequences in
$R$ which converge to $f$, then for all $i\in{\bf N}$, there
exists $n(i)$ such that $f_n-g_l\in m_R^i$ if $n,l\ge n(i)$. Thus
(1)(or (2)) holds for $\{f_n\}$ if and only if (1) (or (2)) holds
for $\{g_l\}$.
\end{pf}

\begin{Definition} Let $R$ be as in the statement of Lemma \ref{LemmaA1}. Let
$$
\begin{array}{l}
p(\hat R)_{\infty}=\\
\{f\in\hat R\mid (2)\text{ of Lemma \ref{LemmaA1} holds for a
Cauchy sequence }\{f_n\}\text{ in }R\text{ which converges to
$f$}\}
\end{array}
$$
\end{Definition}
\begin{Lemma}\label{LemmaA2}
Let $R$ be as in the statement of Lemma \ref{LemmaA1}. Then
\begin{enumerate}
\item $p(\hat R)_{\infty}$ is a prime ideal of $\hat R$ such that
$p(\hat R)_{\infty}\cap R=(0)$. \item There exists a unique
extension $\overline \nu$ of $\nu$ to the quotient field
$\overline K$ of $\hat R/p(\hat R)_{\infty}$ which dominates $\hat
R$. Let $\overline V$ be the valuation ring of $\overline\nu$.
Then $(\overline V,\overline K)$ is an immediate extension of
$(V,K)$. That is, $\Gamma_{\overline V}=\Gamma_{V}$ and
$k(V)=k(\overline V)$.
\end{enumerate}
\end{Lemma}

\begin{pf} The facts that $p(\hat R)_{\infty}$ is prime and $p(\hat R)_{\infty}\cap R=(0)$ are
immediate from Lemma \ref{LemmaA1}.  Suppose that $0\ne f+p(\hat
R)_{\infty}\in\hat R/p(\hat R)_{\infty}$. We can find a Cauchy
sequence $\{f_n\}$ in $R$ such that $\{f_n\}$ satisfies (1) of
Lemma \ref{LemmaA1} and $\{f_n\}$ converges to $f$. Let
$\rho=\nu(f_n)$ for $n>>0$. We necessarily have that $\overline
\nu(g)=\rho$ if $g= f+p(\hat R)_{\infty}$.
\end{pf}

By a classical abuse of notation, we will say that $\nu(f)=\infty$ if $f\in p(\hat R)_{\infty}$.

\section{Simultaneous resolution of $p_{\infty}$}\label{Section4}

\begin{Definition}\label{DefinitionFE2} Suppose that $R$ is a normal local ring which is essentially of
finite type over a field $k$ of characteristic zero, with quotient field $K$. A
normal uniformizing transformation sequence (NUTS) is a sequence of ring homomorphisms
\begin{equation}\label{eqFE1}
\begin{array}{llllllllllllllll}
R&\rightarrow&\overline T_0''&\rightarrow&\overline T_0\\
&&\downarrow&&&\searrow\\
&&\overline T_1'&\rightarrow&\overline T_1''&\rightarrow&\overline T_1\\
&&&&\downarrow&&&\searrow\\
&&&&\overline T_2'&\rightarrow&\overline T_2''&\rightarrow&\overline T_2\\
&&&&&&\downarrow&&&\searrow\\
&&&&&&\vdots&&&\vdots\\
&&&&&&&&&&&\downarrow\\
&&&&&&&&&&&\overline T_n'&\rightarrow&\overline T_n''&\rightarrow&\overline T_n
\end{array}
\end{equation}
such that $\overline T_0=\hat R$, the completion of $R$ with respect to its maximal ideal,
and for all $i$, $\overline T_i$ is the completion with respect to its maximal ideal of
a birational extension $\overline T_i'$ of $\overline T_{i-1}''$. For all $i$, $\overline T_i'$ is a normal local ring, $\overline T_i''$
is a normal local ring, essentially of finite type over $\overline T_i'$ with quotient field
$K_i$ such that $\overline T_i'\subset \overline T_i''\subset \overline T_i$ and $K_0$ is
a finite extension of $K$, $K_{i+1}$ is a finite extension of $K_i$ for all $i\ge 0$.
\end{Definition}

Definition \ref{DefinitionFE2} is the extension of the definition of a UTS in
Chapter 3 of \cite{Dale2} to normal local rings.

To simplify notation, we will often denote the NUTS (\ref{eqFE1}) by $(R,\overline T_n'',\overline T_n)$
or by
$$
R\rightarrow \overline T_0\rightarrow \overline T_1\rightarrow \overline T_n.
$$
We will denote the NUTS consisting of the maps
$$
\begin{array}{lllllll}
\overline T_{n-1}'&\rightarrow&\overline T_{n-1}''&\rightarrow&\overline T_{n-1}\\
&&\downarrow&&&\searrow\\
&&\overline T_n'&\rightarrow&\overline T_n''&\rightarrow&\overline T_n
\end{array}
$$
by $\overline T_{n-1}\rightarrow \overline T_n$.

Suppose that $\nu$ is a  rank 1 $k$-valuation of $K$, and $R$ is
dominated by $\nu$. Suppose that $\overline \nu$ is an extension
of $\nu$ to the quotient field of $\overline T_n$ which dominates
$\overline T_n$. Then we will say that $\overline T_0\rightarrow
\overline T_n$ is a NUTS along $\nu$. When there is no danger of
confusion, we will denote $\overline \nu$ by $\nu$.

We define
$$
p(\overline T_n)_{\infty}=\{f\in \overline T_n\mid \overline
\nu(f)=\infty\},
$$
$$
\lambda(\overline T_n'')=\lambda(\overline T_n)=\text{ dim }\overline T_n/p(\overline T_n)_{\infty}.
$$

We define
$$
\lambda(R)=\lambda(\overline T_0)=\text{ dim }\hat R /p(\hat
R)_{\infty}.
$$
\begin{Lemma}\label{LemmaFE3} Suppose that
$$
\overline T\rightarrow \overline T(1)\rightarrow\cdots\rightarrow \overline T(t)
$$
is a  NUTS along $\nu$. Then $\lambda(\overline T)\ge \lambda(\overline T(t))$.
\end{Lemma}

Lemma \ref{LemmaFE3} is the generalization of Lemma 6.3 of \cite{Dale3} to a NUTS.
The proof is the same.

Let $V$ be the valuation ring (in $K$) of $\nu$ and let
$$
\lambda_V=\text{ min }\{\lambda(R)\mid R \text{ is a normal algebraic local ring of $K$ which
is dominated by $V$}\}.
$$

\begin{Theorem}\label{TheoremFE4} There exists an algebraic regular local ring $\overline R$
of $K$ such that if $(R_1,\overline T''(t),\overline T(t))$ is a
NUTS along $\nu$ with  $R_1$ an algebraic normal local ring of $K$
such that $R_1$ dominates $\overline R $, then
$$
\lambda(\overline T(t))=\lambda_V.
$$
\end{Theorem}

\begin{pf} Suppose that $R$ is an algebraic normal local ring such that $\lambda(R)=\lambda_V$.
Let $\overline R$ be an algebraic regular local ring of $K$ such
that $\overline R$ dominates $R$ and $V$ dominates $\overline R$.
$\lambda(\overline R)=\lambda_V$ by Lemma \ref{LemmaFE3}. Suppose
that $(R_1,\overline T''(t),\overline T(t))$ is a NUTS along $\nu$
with,
 $R_1$ an algebraic normal local ring of $K$ such that
$R_1$ dominates $\overline R $. We have $\lambda(R_1)=\lambda_V$
and $\lambda(\overline T(t))\le\lambda_V$ by Lemma \ref{LemmaFE3}.
Let $L$ be the quotient field of $\overline T''(t)$. $L$ is a
finite extension of $K$. By Theorem 4.2 \cite{Dale} there exists
an  algebraic normal local ring $R_3$ of $K$  and an algebraic
regular local ring $T_2$ of $L$ such that our extension of $V$ to
$L$ dominates $T_2$, $T_2$ dominates  $\overline T''(t)$, $T_2$
dominates $R_3$ and $R_3$ dominates $R_1$, with the property that
$\hat T_2$ is finite over $\hat R_3$. Since $p(\hat
T_2)_{\infty}\cap \hat R_3=p(\hat R_3)_{\infty}$, we have that
$\lambda(R_3)=\lambda(T_2)$. By Lemma \ref{LemmaFE3},
$$
\lambda_V=\lambda(R_1)\ge \lambda(\overline T''(t))\ge \lambda(T_2)=\lambda(R_3)=\lambda_V.
$$
\end{pf}

\begin{Theorem}\label{TheoremFE5}  Suppose that $K^*$ is a finite field extension of $K$
and $\nu^*$ is an extension of $\nu$ to $K^*$.
Then $\lambda_{V^*}=\lambda_V$ and there exists an algebraic  regular local ring $\tilde R$ of $K$ such that the
conclusions of Theorem \ref{TheoremFE4} hold with $\overline R=\tilde R$ and if $S$
is an algebraic normal local ring of $K^*$ which is dominated by $V^*$ and dominates $\tilde R$,\and
$(S, \overline T''(t),\overline T(t))$ is a NUTS along $\nu^*$, then
$\lambda(\overline T(t))=\lambda_V$.
\end{Theorem}

\begin{pf} Let $\overline R$ be the regular local ring of the conclusions of Theorem \ref{TheoremFE4}. Let $S_1$ be a normal algebraic local ring of $K^*$ such that $V^*$ dominates
$S_1$ and $\lambda(S_1)=\lambda_{V^*}$. Let $S_2$ be an algebraic regular local ring such that
$V^*$ dominates $S_2$, $S_2$ dominates $S_1$ and $S_2$ dominates $\overline R$. By Lemma \ref{LemmaFE3}, $\lambda(S_2)
=\lambda_{V^*}$. By Theorem 4.2 \cite{Dale} there exists an  algebraic normal local ring $R_1$
of $K$  and an algebraic regular local ring $S_3$ of $K^*$ such that $V^*$ dominates $S_3$, $S_3$ dominates  $S_2$, $S_3$ dominates
$R_1$ and $R_1$ dominates $\overline R$, with the property that $\hat S_3$ is finite over $\hat R_1$.
Since $p(\hat S_3)_{\infty}\cap \hat R_1=p(\hat R_1)_{\infty}$, $\lambda(R_1)=\lambda(S_3)$.
By Theorem \ref{TheoremFE4} and Lemma \ref{LemmaFE3},
$$
\lambda_V=\lambda(R_1)=\lambda(S_3)=\lambda(S_2)=\lambda_{V^*}.
$$
By Lemma 5.3 \cite{Dale} there exists an algebraic regular local
ring $\tilde R$ of $K$ such that $V$ dominates $\tilde R$, $\tilde
R$ dominates $\overline R$ and if $S$ is an algebraic normal local
ring of $K^*$ which is dominated by $V^*$ and which contains
$\tilde R$, then $S$ dominates $S_2$.

By Theorem \ref{TheoremFE4} applied to $S\subset K^*$ which
dominates $S_2$, the conclusions of Theorem \ref{TheoremFE5} hold.
\end{pf}

We now state a generalization of Theorem 5.1 \cite{Dale2} which resolves
the prime ideal of infinite value terms.

\begin{Theorem}\label{TheoremFE6} Let $k$ be a field of characteristic zero, $K$ an algebraic
function field, $K^*$ a finite algebraic extension of $K$, $\nu^*$ a $k$-valuation of $K^*$,
$\nu=\nu^*\mid K$,
such that $\text{rank }\nu=1$, $\text{rat rank }\nu =s$  and
$$
s\le \lambda=\lambda_V\le n,
$$
where $V$ is the valuation ring of $\nu$ and $n=\text{trdeg}_kK-\text{trdeg}_kV/m_V$.
Suppose that $S^*$ is an algebraic local ring with quotient field $K^*$ which is
dominated by $\nu^*$ and $R^*$ is an algebraic local ring with quotient field $K$
which is dominated by $S^*$. Let $V^*$ be the valuation ring of $\nu^*$.
Then there exists a commutative diagram
$$
\begin{array}{llll}
R_0&\rightarrow&S&\subset V^*\\
\uparrow&&\uparrow\\
R^*&\rightarrow &S^*
\end{array}
$$
where $S^*\rightarrow S$ and $R^*\rightarrow R_0$ are sequences of monoidal transforms
along $\nu^*$
such that $R_0$ has regular parameters $(x_1,\ldots, x_n)$ and $S$ has regular parameters
$(y_1,\ldots,y_n)$  such that there are units $\delta_1,\ldots,\delta_s\in S$ and a
$s\times s$ matrix $A=(a_{ij})$ of natural numbers such that $\text{det}(A)\ne 0$,
\begin{equation}\label{eqFE18}
\begin{array}{ll}
x_1&=y_1^{a_{11}}\cdots y_s^{a_{1s}}\delta_1\\
&\vdots\\
x_s&=y_1^{a_{s1}}\cdots y_s^{a_{ss}}\delta_s\\
x_{s+1}&=y_{s+1}\\
&\vdots\\
x_n&=y_n
\end{array}
\end{equation}
and $\{\nu(x_1),\ldots,\nu(x_s)\}$, $\{\nu(y_1),\ldots,\nu(y_s)\}$
are rational bases of $\Gamma_{\nu}\otimes{\bf Q}=\Gamma_{\nu^*}\otimes
{\bf Q}$. Furthermore,
$$
p(\hat R_0)_{\infty}=(g_{1},\ldots,g_{n-\lambda})
$$
with
$$
g_{i}\equiv x_{s+i}\text{ mod }m(\hat R_0)^2
$$
for $1\le i\le n-\lambda$, and
$$
p(\hat S)_{\infty}=p(\hat R_0)_{\infty}\hat S
$$
are regular primes.
\end{Theorem}

\begin{Remark}\label{RemarkFE19} Suppose that in the hypothesis of
Theorem \ref{TheoremFE6} we further assume that $R^*\rightarrow S^*$ is
such that $R^*$ and $S^*$ have regular parameters $(x_1^*,\ldots,
x_n^*)$, $(y_1^*,\ldots,y_n^*)$ satisfying $(\ref{eqFE18})$ and
such that $\{\nu(x_1^*),\ldots,\nu(x_n^*)\}$,
$\{\nu(y_1^*),\ldots, \nu(y_n^*)\}$ are rational bases of
$\Gamma_{\nu}\otimes{\bf Q}=\Gamma_{\nu^*}\otimes{\bf Q}$. Then in
the conclusions of Theorem \ref{TheoremFE6} we further have that
there exist $b_j(i)\in{\bf N}$, units $\alpha_i\in R_0$ such that
$$
x_i^*=x_1^{b_1(j)}\cdots x_s^{b_s(j)}\alpha_i
$$
for $1\le i\le s$ and there exist $c_j(i)\in{\bf N}$, units
$\beta_i\in S$ such that
$$
y_i^*=y_1^{c_1(j)}\cdots y_s^{c_s(j)}\beta_i
$$
for $1\le i\le s$.

This follows since all transformations in the proof of Theorem
\ref{TheoremFE6} are ``CUTS in the first $n$ variables'' $($page
$49$
 \cite{Dale2}$)$.
\end{Remark}

\begin{pf}(of Theorem \ref{TheoremFE6}). Let $\tilde R$ be the
regular local ring of the conclusions of Theorem
\ref{TheoremFE5}. We first construct a commutative diagram
$$
\begin{array}{llll}
R_1&\rightarrow&S_1&\subset V^*\\
\uparrow&&\uparrow\\
R^*&\rightarrow& S^*
\end{array}
$$
such that the conclusions of Theorem 5.1 \cite{Dale2} hold, and
$R_1$ dominates $\tilde R$. Let  $R=R_1$, $\overline T''=R_1$ and
$\overline T=\hat R_1$. We will now show that we can construct a
CUTS $\overline T\rightarrow \overline T(t)$ along $\nu$, which is
in the first $n$ variables (with the notation of Theorem 4.7 of
\cite{Dale2}), such that $p(\overline T(t))_{\infty}$ has the form
of (53) of page 49 of \cite{Dale2},
\begin{equation}\label{eqFE6}
p(\overline T(t))_{\infty} =(\overline
z_{r(1)}(t)-Q_{r(1)}(\overline z_1(t),\ldots,\overline
z_{r(1)-1}(t)), \ldots, \overline
z_{r(n-\lambda)}(t)-Q_{r(n-\lambda)}(\overline z_1(t),\ldots,\overline
z_{r(n-\lambda)-1}(t)))
\end{equation}
with $s<r(1)<r(2)<\cdots<r(n-\lambda)\le n$ and such that for $1\le
i\le n-\lambda$,
$$
Q_{r(i)}=\overline z_1(t)^{a_1(r(i))}\cdots \overline z_{s}(t)^{a_s(r(i))}u_{r(i)}
$$
where $u_{r(i)}$ is a unit series in $\overline
z_1(t),\ldots,\overline z_{r(i)}(t)$ with coefficients in
$k(c_0,\ldots,c_t)$ (with the notation of (53) of page 49 of
\cite{Dale2}).

The construction of $\overline T\rightarrow \overline T(t)$ follows from the proof of (53) of \cite{Dale2}, with the insertion of the following
at the bottom of page 54. ``Since $\nu(Q_{\overline m}) =\nu(\overline z_{\overline m}(t))<\infty$
we can perform by (54) \cite{Dale2} a UTS in the first $\overline m-1$ variables to get $Q_{\overline m}=\overline z_1(t')^{a_1}\cdots\overline z_s(t')^{a_s}u_{\overline m}$
where $u_{\overline m}\in k(c_0,\ldots,c_{t'})[[\overline z_1(t'),\ldots,\overline z_{\overline m-1}(t')]]$ is a unit series''.

Set $S=S_1$, $\overline U''=S_1$, $\overline U=\hat S$. We can now construct a CUTS $\overline U\rightarrow \overline U(t')$ so that $(R_1,\overline T''(t'),\overline T(t'))$ and
$(S_1,\overline U''(t'),\overline U(t'))$ is a CUTS along $\nu^*$, by Lemma 4.3 and Lemma 4.4
\cite{Dale2}.

Set $\tilde g_j(\overline z_1(t'),\ldots,\overline z_{r(j)}(t'))=\overline z_{r(j)}(t')-
Q_{r(j)}(\overline z_1(t'),\ldots,\overline z_{r(j)-1}(t'))$ for $1\le j\le n-\lambda$.

We will now show that the strict transform of $p(\hat R_1)_{\infty}$ in $\overline T(t')$
is $p(\overline T(t'))_{\infty}$. It suffices to show that the strict transform of
$p(\hat R_1)_{\infty}$ in $\overline T(1)$ is $p(\overline T(1))_{\infty}$. Then the result
follows by induction on $t'$. Let $p=p(\hat R_1)_{\infty}$.  There exists an ideal $I$
in $R_1$, $f\in I$, and a maximal ideal $n$ in $\hat R_1[\frac{I}{f}]$ such that
$\overline T(1)=\widehat{R_1[\frac{I}{f}]_n}$. Let
$$
\overline p=\cup_{j=1}^{\infty}\left(p\hat
R_1[\frac{I}{f}]_n:I^j\hat R_1[\frac{I}{f}]_n\right)
$$
be the strict transform of $p$ in $\hat R_1[\frac{I}{f}]_n$.
$\overline p\ne \hat R_1[\frac{I}{f}]_n$ since the strict
transform in $\hat R_1[\frac{I}{f}]_n$ of an element of infinite
value must have infinite value. $\overline p$ is a prime ideal in
$\hat R_1[\frac{I}{f}]_n$, and
$$
\hat R_1/p\rightarrow \hat R_1[\frac{I}{f}]_n/\overline p
$$
is birational  (Section 0.2 \cite{Hir}, Corollary II.7.15
\cite{Ha}). Thus $\text{dim }\overline T(1)/\overline p\overline
T(1)=\lambda$. $\overline p\overline T(1)$ is a prime contained in
$p(\overline T(1))_{\infty}$ and $\text{dim }\overline
T(1)/p(\overline T(1))_{\infty}=\lambda$ by Theorem
\ref{TheoremFE5} (since $R_1$ contains $\tilde R$). Thus
$\overline p\overline T(1)=p(\overline T(1))_{\infty}$.

Hence there exist $f_1,\ldots, f_{n-\lambda}\in p(\hat
R_1)_{\infty}$, $c_1(i),\ldots,c_s(i)\in{\bf N}$ for $1\le i\le
n-\lambda$ and $b_{ij}\in \overline T(t')$ such that
$\text{det}(b_{ij})$ is a unit in $\overline T(t')$,
$f_i=M_i(\sum_{j=1}^{n-\lambda} b_{ij}\tilde g_j)$ where
$M_i=\overline z_1(t')^{c_1(i)}\cdots\overline z_s(t')^{c_s(i)}$
for $1\le i\le n-\lambda$ and
$$
p(\overline T(t'))_{\infty}=\left(\frac{f_1}{M_1},\ldots,\frac{f_{n-\lambda}}{M_{n-\lambda}} \right).
$$
Let $m$ be a positive integer such that
\begin{equation}\label{eqFE10}
m> \frac{\left(\text{max}_{1\le j\le n-\lambda}\nu(\overline z_{r(j)}(t'))\right)
+\left(\text{max}_{1\le i\le n-\lambda}\nu(M_i)\right)}{\nu(m(\overline T(t')))}.
\end{equation}

By Theorem 4.8 \cite {Dale2} (with $l=n$) there exists a CRUTS along $\nu$,
$(R_1,R_1,T(t'))$ and $(S_1,S_1,U(t'))$ with associated MTSs
$$
\begin{array}{lll}
S&\rightarrow &S(t')\\
\uparrow&&\uparrow\\
R&\rightarrow&R(t')
\end{array}
$$
such that (with the notation of Theorem 4.8 \cite{Dale2})
$$
f_i=\overline M_i(\sum_{j=1}^{n-\lambda} b_{ij}(\tilde{\overline x}_1(t'),\ldots,\tilde{\overline x}_n(t'))
\tilde g_j(\tilde{\overline x}_1(t'),\ldots,\tilde{\overline x}_{r(j)}(t'))+h_i\in p(\hat R(t'))_{\infty}
$$
with $\overline M_i=\tilde{\overline
x}_1(t')^{c_1(i)}\cdots\tilde{\overline x}(t')^{c_s(i)}$,
 $h_i\in m(T(t'))^m$ (where $m$ is the integer of (\ref{eqFE10}))
for $1\le i\le n-\lambda$, and such that (by (A3) of page 83 of \cite{Dale2})
$\nu(\overline z_i(t'))=\nu(\tilde{\overline x}_i(t'))$ for $1\le i\le n$, and
$$
\begin{array}{ll}
\nu(\overline z_{r(i)}(t'))&=\nu(Q_{r(i)}(\overline z_1(t'),\ldots,\overline z_{r(i)-1}(t')))\\
&=\nu(Q_{r(i)}(\tilde{\overline x}_1(t'),\ldots,\tilde{\overline x}_{r(i)-1}(t')))\\
&=\nu(\tilde{\overline x}_{r(i)}(t'))
\end{array}
$$
for $1\le i\le n-\lambda$.

Now we perform the MTS
$$
\begin{array}{lll}
S(t')&\rightarrow & S(t'')\\
\uparrow&&\uparrow\\
R(t')&\rightarrow& R(t'')
\end{array}
$$
of the proof of Theorem  4.9 \cite{Dale2} (with $l=n$). Because of the form of the $\tilde g_j$, we have for $1\le j\le n-\lambda$,
$$
\tilde g_j=\overline x_1(t'')^{b_1(j)}\cdots \overline x_s(t'')^{b_s(j)}\overline g_j
$$
where
$$
\overline g_j\equiv \tilde d_j\overline x_{r(j)}(t'')\text{ mod }(\overline x_1(t''),\ldots,\overline x_{r(j)-1}(t''))\hat R(t'')
+m(\hat R(t''))^2
$$
for some nonzero $\tilde d_j\in R(t'')/m(R(t''))$, where
\begin{equation}\label{eqFE11}
\nu(\overline x_1(t'')^{b_1(j)}\cdots \overline x_s(t'')^{b_s(j)})=\nu(\tilde{\overline x}_{r(j)}(t'))=\nu(\overline z_{r(j)}(t')).
\end{equation}

We further have
$$
\overline h_i=\frac{h_i}{\overline M_i}=\overline x_1(t'')^{e_1(i)}\cdots \overline x_s(t'')^{e_s(i)}\tilde h_i
$$
for $1\le i\le n-\lambda$ with $\tilde h_i\in \hat R(t')$ and
\begin{equation}\label{eqFE12}
\nu(\overline x_1(t'')^{e_1(i)}\cdots \overline x_s(t'')^{e_s(i)})>\text{ max }\nu(\overline z_{r(j)}(t))
\end{equation}
for $1\le j\le n-\lambda$ by (\ref{eqFE10}). By Lemma 4.2
\cite{Dale2}, (\ref{eqFE11}) and (\ref{eqFE12}), we can further
choose the final CUTS of type (M1) (on the top of page 89 of
\cite{Dale2}) so that
$$
e_j(i)>\text{ max}_{1\le \alpha\le n-\lambda}b_j(\alpha)
$$
for $1\le j\le s$ and $1\le i\le n-\lambda$.

Let $A=(b_{ij})^{-1}$, a matrix with coefficients in $\hat R(t'')$.
$$
A\left(\begin{array}{l}
\frac{f_1}{\overline M_1}\\ \vdots\\  \frac{f_{n-\lambda}}{\overline M_{n-\lambda}}\end{array}\right)
=
\left(\begin{array}{l}\tilde g_1+d_1\\ \vdots\\ \tilde g_{n-\lambda}+d_{n-\lambda}\end{array}
\right)
$$
where for $1\le i\le n-\lambda$,
$$
d_i=\overline x_1(t'')^{b_1(i)+1}\cdots \overline x_s(t'')^{b_s(i)+1}\overline d_i
$$
for some $\overline d_i\in \hat R(t'')$. Thus
$$
g_i= \overline g_i+\overline d_i\in p(\hat R(t''))_{\infty}
$$
 and
$$
g_i\equiv \tilde d_i\overline x_{r(i)}(t'')\text{ mod }(\overline x_1(t''),\ldots,
\overline x_{r(i)-1}(t''))\hat R(t'')+m(\hat R(t''))^2.
$$
Thus $(g_1,\ldots, g_{n-\lambda})$ is a complete intersection and a
regular prime ideal in $\hat R(t'')$. Since $\lambda(\hat
R(t''))=\lambda$ (by Theorem \ref{TheoremFE4}), we have that
$(g_1,\ldots,g_{n-\lambda})$ is a basis of $p(\hat
R(t''))_{\infty}$. Since $(g_1,\ldots,g_{n-\lambda})\hat S(t'')$ is
a prime ideal and $\lambda(S(t''))=\lambda$ (by Theorem
\ref{TheoremFE5}), it follows that $p(\hat R(t''))_{\infty}\hat
S(t'') =p(\hat S(t''))_{\infty}$.

We can now make a change of variables in the regular parameters $(x_1(t''),\ldots,x_n(t''))$
and $(y_1(t''),\ldots,y_n(t''))$ of the proof of Theorem 4.9 \cite{Dale2} to get the desired
forms of the $g_i$.

$\{\nu(x_1),\ldots,\nu(x_s)\}$ and $\{\nu^*(y_1),\ldots,\nu^*(y_s)\}$ are rational bases of
$\Gamma^*\otimes{\bf Q}$ by the construction of the sequence $R^*\rightarrow R_0$ and
$S^*\rightarrow S$.
\end{pf}

\section{Rank 1 valuations}\label{Section5}
Let notations be as in Section \ref{Section2}. Further assume that
$V$ has rank 1. Consider our directed set $\{R_i\mid i\in I\}$
satisfying (a), (b), (c) and (d). For $i\in I$, we define
$$
p_i=p(\hat R_i)_{\infty}=\{f\in \hat R_i\mid \nu(f)=\infty\}.
$$
For $i<j$, the natural inclusions $\hat R_i\rightarrow \hat R_j$
induce inclusions $\hat R_i/p_i\rightarrow \hat R_j/p_j$. Thus
$\{\hat R_i/p_i\mid i\in I\}$ is a directed system, and we have a
local domain
$$
\overline T=\lim_{\rightarrow}\hat R_i/p_i
=\cup\hat R_i/p_i.
$$
Let $\overline K_{\infty}$ be the quotient field of $\overline T$.

\begin{Lemma}\label{lemma} Suppose that $V$ has rank 1. Then
the ring $\overline T=\cup_{i\in I}\hat {R_i}/p_i$ does not depend
on the directed system of rings $\{R_i\mid i\in I\}$ satisfying
$(a)$, $(b)$, $(c)$ and $(d)$.
\end{Lemma}

\begin{pf}
The proof is essentially the same as that of Lemma \ref{LemmaA5}. We must observe that
the inclusion $\hat R_i\rightarrow \hat S_j$ of the proof of Lemma \ref{LemmaA5}
induces a natural inclusion $\hat R_i/p_i\rightarrow \hat S_j/q_j$,
where $q_j=\{f\in \hat S_j\mid \nu(f)=\infty\}$.
\end{pf}

Theorem \ref{TheoremA3} is a generalization of Zariski's local uniformization theorem
 \cite{Z}.  Our proof is an extension in rank 1 of the proof for general rank in
\cite[6.2]{Dale}. We incorporate the conclusions of Theorem
\ref{TheoremFE6} which resolves the prime ideal of infinite value
terms.

\begin{Theorem}\label{TheoremA3}
Let $k$ be a field of characteristic zero, $K$ an algebraic
function field over $k$, and let $\nu$ be a rank 1 $k$-valuation
of $K$, of rational rank $s$, with valuation ring $V$. Let
$$
n=\text{trdeg}_kK-\text{trdeg}_kV/m_V,
$$
$\lambda=\lambda_V$ $($defined before Theorem \ref{TheoremFE4}$)$.

 Then there exists a partially ordered set $I$ and
algebraic regular local rings $\{R_i\mid i\in I\}$ with quotient
field $K$ which are dominated by $V$ such that
$$
V=\lim_{\rightarrow}R_j=\cup_{j\in I} R_j
$$
and $R_j$ has regular parameters $(x_1(j),\ldots,x_n(j))$ such
that
\begin{enumerate}
\item
$$
\{\nu(x_1(j)),\ldots,\nu(x_{s}(j)\}
$$
is a rational basis of $\Gamma\otimes{\bold Q}$.
\item
If $j<k\in I$ then there are relations
\begin{equation}\label{eqV1}
x_i(j)=\prod_{l=1}^sx_l(k)^{d_{il}}\delta_i
\end{equation}
for $1\le i\le s$ where $\delta_i\in R_k$ are units. The $s\times
s$ matrix $D(j,k)=(d_{il})$ of $(\ref{eqV1})$ has nonzero
determinant. \item The prime ideal
$$
p_j=p(\hat R_j)_{\infty}=\{f\in \hat R_j\mid \nu(f)=\infty\}
=(g_{1}(j),\ldots, g_{n-\lambda}(j))
$$
with
$$
g_{i}(j)\equiv x_{s+i}(j)\text{ mod } m(\hat R_j)^2.
$$
In particular, $p_j$
is a regular prime.
\item
For $j\in I$, let $\Lambda_j$ be the free ${\bold Z}$-module
$\Lambda_j=\sum_{i=1}^s\nu(x_i(j)){\bold Z}$. Then
$$
\Gamma = \lim_{\rightarrow}\Lambda_j=\cup_{j\in I}\Lambda_j.
$$
\end{enumerate}
\end{Theorem}

\begin{pf} Let $R^*$ be an algebraic regular local ring such that $V$ dominates $R^*$.
By
 Theorem \ref{TheoremFE6} (with $K=K^*$ and $R^*=S^*$), there exists a sequence of monoidal
transforms $R^*\rightarrow R_0$ along $V$ such that (1) and (3) of
this theorem hold on $R_0$.

Suppose that $m$ is a positive integer and $f=(f_1,\ldots,f_m)\in
V^m$. We will construct a sequence of monoidal transforms
$R_0\rightarrow R_f$ along $V$ such that $f_1,\ldots, f_m\in R_f$,
(1) and (3) of this theorem hold for $R_f$ and (2) of this theorem
holds for $R_0\rightarrow R_f$. We will further have
$\nu(f_1),\ldots, \nu(f_m)\in \Lambda_f$.

 By Theorem 4.9 \cite{Dale}
with the $R^*$, $S^*$ of the statement of Theorem 4.9 set as
$R^*=S^*=R_0$, and $v_i=x_i(0)$ if  $1\le i\le s$, and
$v_{s+1}=f_1,\ldots, v_{s+m}=f_m$,
 there exists a sequence of
monoidal transforms $R_0\rightarrow R_1$ along $V$ such that (1)
of this theorem holds for $R_1$, (2) of this theorem holds for
$R_0\rightarrow R_1$, $f_1,\ldots,f_m\in R_1$ and
$\nu(x_1),\ldots,\nu(x_s)\in\Lambda_1$.  By Theorem
\ref{TheoremFE6} (with $K=K^*$, $R^*=S^*$) and Remark
\ref{RemarkFE19}, there exists a sequence of monoidal transforms
$R_1\rightarrow R_f$ along $V$ such that (1), (2) and (3) of this
theorem hold for $R_0\rightarrow R_f$ and
$\nu(f_1),\ldots,\nu(f_m)\in\Lambda_f$. We have
$\text{det}D(0,f)\ne 0$ since $\{\nu(x_i(0))\mid 1\le i\le s\}$
and $\{\nu(x_i(f))\mid 1\le i\le s\}$ are two bases of
$\Gamma\otimes {\bold Q}$.

Let $I=\sqcup_{m\in {\bold N}_+}V^m$ be the disjoint union. For
$f\in I$ we construct $R_f$ as above. If $f=0$ we let $R_0$ be the
$R_0$ constructed above. Define a partial order on $I$ by $f\le g$
if $R_f\subset R_g$.

Suppose that $R_{\alpha}\subset R_{\beta}$. We have $R_0\subset
R_{\alpha}\subset R_{\beta}$.
$$
x_i(0)=\prod_{j=1}^sx_j(\alpha)^{c_{ij}}\delta_i
$$
for $1\le i\le s$ with $\delta_i$ a unit in $R_{\alpha}$ and
$$
x_i(0)=\prod_{j=1}^sx_j(\beta)^{d_{ij}}\epsilon_i
$$
for $1\le i\le s$ with $\epsilon_i$ a unit in $R_{\beta}$. Thus in
 $R_{\beta}$ there are factorizations
$$
x_i(\alpha)=\prod_{j=1}^sx_j(\beta)^{e_{ij}}\lambda_i
$$
for $1\le i\le s$ and $\lambda_i$ a unit in $R_{\beta}$. We have
$\text{det}(D(\alpha,\beta))\ne 0$ since (1) holds for
$R_{\alpha}$ and $R_{\beta}$.  Thus (2) holds for $R_{\alpha}
\rightarrow R_{\beta}$. To show that $V=\lim_{\rightarrow}R_j$, we
must verify that $I$ is a directed set. That is, for
$\alpha,\beta\in I$, there exists $\gamma\in I$ such that
$R_{\alpha}\subset R_{\gamma}$ and $R_{\beta}\subset R_{\gamma}$.

There exists $f_1,\ldots,f_t\in V$ such that if
$A=k[f_1,\ldots,f_t]$, $m=A\cap m_{V}$ then $R_{\alpha}=A_m$.
There exists $g_1,\ldots, g_n\in V$ such that if
$B=k[g_1,\ldots,g_n]$, $n=B\cap m_{V}$ then $R_{\beta}=B_n$. Set
$\gamma=(f_1,\ldots,f_t,g_1,\ldots,g_n)$. By construction, $A,
B\subset R_{\gamma}$. Since $m_V\cap R_{\gamma}=m_{\gamma}$ is the
maximal ideal of $R_{\gamma}$, we have $R_{\alpha}$,
$R_{\beta}\subset R_{\gamma}$.

 (4) holds by our construction, since $\nu(f)\in\Lambda_f$ if $f\in V$.
\end{pf}

\medskip

\begin{Theorem}\label{TheoremA6} Suppose that $V$ has rank 1 and $k$ has characteristic zero. Then
 the ring $\overline T=\cup_{i\in I}(\hat {R_i}/p_i)$ is a valuation ring.
\end{Theorem}

\begin{pf} Let $s$ denote the rational rank of $\nu$.
By Lemma \ref{lemma} we can assume that the rings $R_i$'s are as
in Theorem \ref{TheoremA3}. $R_i$ has regular parameters
$x_1(i),\dots, x_n(i)$ and $\hat
{R_i}=R_i/{m}_{R_i}[[x_1(i),\dots, x_n(i)]].$

Let $f\in \hat {R_i}$. We recall that if $\nu(f)<\infty$, then by
Theorem 4.8 and Theorem 4.10 \cite{Dale2} after a MTS
$R_i\rightarrow \overline R(1)$ along $\nu$
$$f= \overline x_1(1)^{d_1}\cdots \overline x_s(1)^{d_s} u( \overline x_1(1),\ldots,
\overline x_n(1))$$ where $ \overline x(1),\ldots, \overline
x_n(1)$ are regular parameters in the ring $\overline R(1)$, $
\nu(\overline x(1)),\ldots, \nu(\overline x_s(1))$ are rationally
independent and $u\in \widehat {\overline R(1)}$ is a unit power
series. Further, there exist units $\alpha_j\in \overline R(1)$
such that
$$
x_j(i)=\overline x_1(1)^{\beta_1(j)}\cdots\overline x_s(1)^{\beta_s(j)}\alpha_j
$$
for $1\le j\le s$.

Let $h\in Q(\overline T)$. We want to show that either $h\in \overline T$, or $1/h \in
\overline T$. So it suffices to show that if $\nu(h)\geq 0$, then $h\in \overline T$.
Write $h=a/b$ where $a\in \cup_{i\in I}\hat {R_i}/p_i$ and $0\neq
b\in \cup_{i\in I}\hat {R_i}/p_i$. Then $a\in \hat {R_j}/p_j$ for
some $j\in I$, and $b\in \hat {R_k}/p_k$ for some $k\in I$. After
a MTS $R_j\rightarrow \overline R(1)$ along $\nu$ we have
$$a= \overline x_1(1)^{d_1}\cdots \overline x_s(1)^{d_s} u(
\overline x_1(1),\ldots, \overline x_n(1))$$ where $ \overline
x_1(1),\ldots, \overline x_n(1)$ are regular parameters in the
ring $ \overline R(1)$ and $u\in \widehat {\overline R(1)}$ is a
unit. Further, there exist units $\alpha_j\in \overline R(1)$ such
that
\begin{equation}\label{eqFE21}
x_j(i)=\overline x_1(1)^{\beta_1(j)}\cdots\overline x_s(1)^{\beta_s(j)}\alpha_j
\end{equation}
for $1\le j\le s$.

After another MTS $R_k\rightarrow \overline R(2)$ along $\nu$ we
have
$$b= \overline x_1(2)^{e_1}\cdots
\overline x_s(2)^{e_s} u'( \overline x_1(2),\ldots, \overline x_n(2))$$ where $
\overline x_1(2),\dots, \overline x_n(2)$ are regular parameters in the ring $
\overline R(2)$ and $u'\in \widehat {\overline R(2)}$ is a unit.
Further, there exist units $\gamma_j\in \overline R(2)$ such that
\begin{equation}\label{eqFE22}
x_j(k)=\overline x_1(2)^{\delta_1(j)}\cdots\overline x_s(2)^{\delta_s(j)}\gamma_j
\end{equation}
for $1\le j\le s$.

We have that $\overline R(1)=k[f_1,\dots,f_m]_{{ m}_V\cap
k[f_1,\dots,f_m]}$ for some $f_1,\dots,f_m \in V$, and $\overline
R(2)=k[g_1,\dots,g_n]_{{ m}_V\cap k[g_1,\dots,g_n]}$ for some
$g_1,\dots,g_n \in V$.

Let $c=(f_1,\dots,f_m,g_1,\dots,g_n)\in V^{m+n}$. Let $R_c$ be
constructed as in the proof of Theorem \ref{TheoremA3}. Then
$R_c\in \cup_{i\in I} R_i$ and $\overline R(1)\subset R_c$,
$\overline R(2)\subset R_c$. The ring $R_c$ has regular parameters
$x_1(c),\dots, x_n(c)$ and by (2) of Theorem \ref{TheoremA3} and
(\ref{eqFE21}) and (\ref{eqFE22}) the ``good form'' of $a$ and $b$
is preserved in $R_c$:
$$a=
x_1(c)^{a_1}\dots x_s(c)^{a_s} \bar u(x_1(c),\dots, x_n(c))$$

$$b=
x_1(c)^{b_1}\dots x_s(c)^{b_s} \bar u'(x_1(c),\dots, x_n(c))$$

where $\bar u$ and $\bar u'$ are units in $\hat {R_c}$.

Let $g=x_1(c)^{a_1}\dots x_s(c)^{a_s}/x_1(c)^{b_1}\dots
x_s(c)^{b_s}$. Since $g\in K=Q(R_c)$ and $\nu(g)=\nu(h)\ge 0$, we
have that $g\in V$ and  $g\in R_g$, which is in the directed
system $I$.  We have $g=x_1(g)^{f_1}\dots
x_s(g)^{f_s}w(x_1(g),\dots x_n(g))$ with $f_i\geq 0$ for every $i$
and $w$ a unit in $R_g$. There exists $\gamma\in I$ such that
$R_c\subset R_{\gamma}$ and $R_g\subset R_{\gamma}$. Then
$$h=x_1(\gamma)^{t_1}\dots x_s(\gamma)^{t_s}\bar w(x_1(\gamma),\dots,
x_n(\gamma)), $$ with $t_i\geq 0$ for every $i$ and $\bar w$ a
unit in $\widehat R_{\gamma}$. Hence $h\in \widehat R_{\gamma}$
and so
 $h\in \overline T$.
\end{pf}

\medskip

\begin{Theorem}\label{TheoremA4} Suppose that $V$ has rank 1 and $k$ has characteristic zero. Then
\begin{enumerate}
\item  $(\overline T, Q(\overline T))$ is an Henselian immediate
extension of $(V,K)$. \item There exists a directed system of
regular local rings $\{R_i\}$ satisfying $(a)$, $(b)$, $(c)$, and
$(d)$ such that each $\hat R_i/p_i$ is a regular local ring, and
$\overline T =\cup \hat R_i/p_i$.
\end{enumerate}
\end{Theorem}

\begin{pf} Let $\overline K_i$ be the quotient field of $\hat R_i/p_i$.
Then for all $i\in I$, $(\overline T\cap \overline K_i,\overline K_i)$ is an immediate extension of
$(V,K)$ by Lemma \ref{LemmaA2}. Thus $(\overline T,Q(\overline T))$
is an immediate extension of $(V,K)$.

By an extension of Proposition \ref{PropV1}, $\overline T$ is Henselian.

Statement (2) follows from the construction of Theorem
\ref{TheoremA3}.
\end{pf}

Suppose that $W$ is a rank 1 valuation ring, with valuation
$\omega$. Let $\phi(x)=e^{-\omega(x)}$ for $x\in L=Q(W)$. A
sequence $(x_i)_{i\in{\bf N}}$ of elements of $L$ is $\phi$-Cauchy
if given $\epsilon >0$, there exists $n_0$ such that
$\phi(x_n-x_m)<\epsilon$ for all $m,n\ge n_0$ (Section 2,
\cite{E}). $W$ is said to be complete if all $\phi$-Cauchy
sequences $(x_i)$ converge to an $x\in L$.

\begin{Example}\label{ExampleFE25}  Even if $V$ has rank 1 and $T=\overline T$ is a valuation ring,
$($which necessarily has residue field $k(T)=k(V)$ and value group
$\Gamma_V=\Gamma_T$$)$, $(T,K_{\infty})$ is not in general
complete. In particular, it is not a maximal immediate extension
$($in the sense of Krull \cite{Kr2} and Kaplansky \cite{K}$)$.
\end{Example}

\begin{pf} Let $K=k(x,y)$ be a rational function field in two variables
over a field $k$ of characteristic zero.
Let $R=k[x,y]_{(x,y)}$. Let $\nu$ be the rank one valuation of $K$
with nondiscrete value group which we can take to be ${\bf Q}$ and residue field $k$ which dominates
$R$ constructed in Example 3, page 102 of \cite{ZS}. Let
$$
R\rightarrow R_1\rightarrow R_2\rightarrow \cdots\rightarrow
R_i\rightarrow\cdots
$$
be the system of regular local rings for $i\in{\bf N}$ of the
construction such that  $\cup R_i=V$ is the valuation ring of
$\nu$.

We will first establish that $T= \cup \hat R_i$ is a valuation
ring with residue field $k(V)$ and value group ${\bf Q}$. For any
fixed $i$, let $q_i=\{f\in R_i\mid\nu(f)=\infty\}$. By Lemma \ref{LemmaA2} $\nu$ extends
uniquely to a valuation of $Q(\hat R_i/q_i)$ which dominates $\hat
R_i/q_i$ and has residue field $k$ and value group ${\bf Q}$. If
$q_i\ne (0)$, then $\hat R_i/q_i$ is a 1 dimensional excellent
local ring, so the only valuation rings of $Q(\hat R_i/q_i)$ which
dominate $\hat R_i/q_i$ are discrete, which is a contradiction.
Thus $q_i=0$. By Theorems \ref{TheoremA6} and \ref{TheoremA4},  $T=\overline T$ is a
valuation ring with value group ${\bf Q}$ and
residue field $k$.

In each regular local ring $\hat R_i$ there is the sequence of all
valuation ideals
$$
\cdots\subset I_n(i)\subset\cdots\subset I_1(i)\subset
I_0(i)=m(\hat R_i)\subset\hat R_i.
$$
Let $p_j(i)=\nu(I_j(i))$. For fixed $j$,
$lim_{j\rightarrow\infty}p_j(i)=\infty$ (c.f. Lemma 2.3
\cite{Dale2}).
 Notice that $\cup_{i\in{\bf
N}}\cup_{j\in{\bf N}}\{p_j(i)\} ={\bf Q}_+$, but arbitrarily large elements of ${\bf Q}_+$
are not in $\cup_{j\in{\bf N}}\{p_j(i)\}$ for a fixed $i$ since $\cup_{j\in{\bf N}}\{p_j(i)\}$
is not discrete.

We can inductively construct for all $i\in{\bf N}$,
$\sigma(i),\lambda(i)\in{\bf N}$ and $a_i\in
I_{\lambda(i)}(\sigma(i))$ such that $i<\sigma(i)$, $\sigma(i-1)<\sigma(i)$,
$\nu(a_i)=p_{\lambda(i)}(\sigma(i))$,
$$
\nu(a_i)>\text{ max
}\{p_{\lambda(0)}(\sigma(0)),\ldots,p_{\lambda(i-1)}(\sigma(i-1)),i\},
$$
 and $p_{\lambda(i)}(\sigma(i))\not\in\cup_{k<\sigma(i)}\cup_{j\in{\bf N}}\{p_j(k)\}$ for every $i$.

For $i\in{\bf N}$, set $\alpha_i=a_1+\cdots+a_i$.  For $i<j$, we
have
$$
\nu(\alpha_j-\alpha_i)=\nu(a_{i+1})>i.
$$
Thus $\{\alpha_i\}$ is a $\phi$-Cauchy sequence.  Suppose that
there exists a limit $\tau\in K_{\infty}$ of $\{\alpha_i\}$. Then
\begin{equation}\label{eq1}
\nu(\tau-\alpha_i)=p_{\lambda(i+1)}(\sigma(i+1))
\end{equation}
for all $i$, by the definition of a limit.  We have $\tau\in T$ so
that $\tau\in \hat R_{\sigma(i)}$ for some $i$.  Thus $\tau-\alpha_i\in\hat
R_{\sigma(i)}$. But by (\ref{eq1}) we have that
$\nu(\tau-\alpha_i)$ is not the value of an element of $\hat
R_{\sigma(i)}$, a contradiction.
\end{pf}
\medskip

\section{Ramification of completions of rank 1 valuation rings}\label{Section6}

Theorem \ref{Theorem8} is a generalization of Theorem 6.3 \cite{Dale}, which resolves
the prime ideal of infinite value terms.

\begin{Theorem}\label{Theorem8}Let $k$ be a field of characteristic zero, $K$ an algebraic function field
over $k$, $K^*$ a finite algebraic extension of $K$, $V^*$ a rank
1 $k$-valuation ring of $K^*$ of rational rank s, $V=V^*\cap K$.
Let
$$
e=[\Gamma^*:\Gamma]
$$
be the ramification index of $V^*$ relative to $V$,
$$
f=[V^*/m_{V^*}:V/m_V]
$$
be the residue degree of $V^*$ relative to $V$, and let $\tau$ be a primitive
element of $V^*/m_{V^*}$ over $V/m_V$. Let
$$
n=\text{trdeg}_kK^*-\text{trdeg}_kV^*/m_{V^*}=\text{trdeg}_kK-\text{trdeg}_kV/m_V,
$$
$\lambda=\lambda_V=\lambda_{V^*}$ $($as shown in Theorem
\ref{TheoremFE5}$)$.

Then there exists a partially ordered set $I$ and algebraic
regular local rings $\{S_i\mid i\in I\}$ with quotient field $K^*$
which are dominated by $V^*$ where $S_j$ has regular parameters
$(y_1(j),\ldots,y_n(j))$ such that
\begin{enumerate}
\item
$$
\{\nu^*(y_{1}(j)),\ldots,\nu^*(y_s(j))\}
$$
is a rational basis of $\Gamma^*\otimes{\bold Q}$.
\item For all $k\in I$ there exist algebraic regular local rings $R_0(k)$
with quotient field $K$ which are dominated by $V$ such that there exist factorizations
$$
R_0(k)\rightarrow R_k\rightarrow S_k
$$
so that there are regular parameters $(x_1(k), .... ,x_n(k))$ in
$R_0(k)$,  units $\delta_1(k),\ldots,\delta_s(k)\in S_k$ and a
$s\times s$ matrix $A(k)=(a_{ij}(k))$ of nonnegative integers such
that  $\text{det}(A(k)) \ne 0$ and
 \begin{equation}\label{eqV3}
\begin{array}{ll}
x_1(k) &= y_1(k)^{a_{11}(k)} \cdots y_s(k)^{a_{1s}(k)}\delta_1(k)\\
&\vdots\\
x_s(k) &= y_1(k)^{a_{s1}(k)} \cdots y_s(k)^{a_{ss}(k)}\delta_s(k)\\
x_{s+1}(k)&= y_{s+1}(k)\\
&\vdots\\
x_n(k)&=y_n(k).
\end{array}
\end{equation}

$R_k$ is a normal local ring with quotient field $K$ $($which is
obtained by a toric blowup of $R_0(k))$ such that $S_k$ is a
localization at a maximal ideal of the integral closure of $R_k$
in $K^*$. The prime ideals
$$
p_j=p(\hat R_j)_{\infty}=\{f\in \hat R_j\mid \nu(f)=\infty\}
=(g_{1}(j),\ldots, g_{n-\lambda}(j))
$$
with
$$
g_{i}(j)\equiv x_{s+i}(j)\text{ mod } m(\hat R_j)^2
$$
and
$$
q_j=p(\hat S_j)_{\infty}=p(\hat R_j)_{\infty}\hat
S_j=(g_{1}(j),\ldots, g_{n-\lambda}(j))
$$
with
$$
g_i(j)\equiv y_{s+i}(j)\text{ mod }m(\hat S_j)^2.
$$

Furthermore, there are
isomorphisms of abelian groups
$$
\Gamma^*/\Gamma\cong {\bold Z}^s/{A(k){\bold Z}^s},
$$
$$
[S_k/m_{S_k}:R_k/m_{R_k}]=f,\,\,\,
\mid\text{det}(A(k))\mid=e,\,\,\, [QF(\hat S_k):QF(\hat R_k)]=ef
$$
and $S_k/m_{S_k}=R_k/m_{R_k}[\tau]$. \item
 Let
$k'$ be an algebraic closure of $V^*/m_{V^*}$.
Suppose that  $j<k\in I$.
\begin{enumerate}
\item There are relations
\begin{equation}\label{eqV2}
y_i(j)=\prod_{c=1}^sy_c(k)^{d_{ic}}\epsilon_i
\end{equation}
where $d_{ic}$ are natural numbers and $\epsilon_i\in S_k$ is a
unit for $1\le i\le s$. Let $D(j,k)$ be the $s\times s$ matrix of
$(\ref{eqV2})$. Then $\text{det}(D(j,k))\ne 0$. \item There exists
a commutative diagram
\begin{equation}\label{eqV4}
\begin{array}{lll}
R_k&\rightarrow&S_k\\
\uparrow&&\uparrow\\
R_j&\rightarrow&S_j
\end{array}
\end{equation}
\item We have actions of $\Gamma^*/\Gamma$ on $\hat
S_j\otimes_{S_j/m_{S_j}}k'$ such that
$$
(\hat S_j\otimes_{S_j/m_{S_j}}k')^{\Gamma^*/\Gamma}\cong \hat
R_j\otimes_{R_j/m_{R_j}}k'
$$
for all $j$, and  this action is compatible with restriction.
 \end{enumerate}
We have an isomorphism
$$
\hat S_j\otimes_{S_j/m_{S_j}}k'\cong k'[[\overline
y_1(j),\ldots,\overline y_n(j)]]
$$
where $\overline y_1(j),\ldots,\overline y_n(j)$ are defined by
$$
x_i(j)=\left\{\begin{array}{ll}
\overline y_1(j)^{a_{i1}(j)}\cdots\overline y_s(j)^{a_{is}(j)}&\text{ if }
1\le i\le s\\
\overline y_i(j)&\text{ if }s+1\le i\le n\end{array}\right.
$$
Let $(b_{\alpha,\beta}(j))=\text{adj} A(j)$ and $\omega$ be a
primitive $e$-th root of unity. The action of $\Gamma^*/\Gamma$ on
$\hat S_j\otimes_{S_j/m_{S_j}}k'\cong k'[[\overline
y_1(j),\ldots,\overline y_n(j)]]$ is defined for
$$
c\in {\bold Z}^s/A(j){\bold Z}^s\cong \Gamma^*/\Gamma
$$
by
$$
\sigma_c(\overline y_{\alpha}(j))=\left\{
\begin{array}{ll}
\omega^{\sum_{\beta=1}^sb_{\alpha,\beta}(j)c_{\beta}}\overline y_{\alpha}(j)&\text{ if }1\le \alpha\le s\\
\overline y_{\alpha}(j)&\text{ if }s+1\le j\le n.\end{array}\right.
$$
\item
$$
V^*=\lim_{\rightarrow}S_j=\cup_{j\in I} S_j
$$
and
$$
V=\lim_{\rightarrow}R_j=\cup_{j\in I} R_j.
$$
For $j\in I$, let $\Lambda_j$ be the free ${\bold Z}$ module
$\Lambda_j=\sum_{i=1}^s\nu(x_i(j)){\bold Z}$, and let
$\Omega_j$ be the free ${\bold Z}$ module
$\Omega_j=\sum_{i=1}^s\nu^*(y_i(j)){\bold Z}$. Then
$$
\Gamma = \lim_{\rightarrow}\Lambda_j=\cup_{j\in I}\Lambda_j
$$
and
$$
\Gamma^* = \lim_{\rightarrow}\Omega_j=\cup_{j\in I}\Omega_j.
$$
\end{enumerate}
\end{Theorem}

\begin{pf}
Suppose that $R'$ is the regular local ring of Theorem 6.1 \cite{Dale},
and $\tilde R$ is the regular local ring of Theorem \ref{TheoremFE5}. By Theorem
6.1 \cite{Dale}, there exists a sequence of local rings
$$
R_0(0)\rightarrow R_0\rightarrow S_0
$$
such that $R'\subset R_0(0)$, $\tilde R\subset R_0(0)$ and the
conclusions of Theorem 6.3 \cite{Dale} and Theorem
\ref{TheoremFE6} hold for this sequence. In particular, (1) and
(2) of the theorem hold for $R_0(0)\rightarrow R_0\rightarrow S_0$
and $p(\hat S_0)_{\infty}$, $p(\hat R_0)_{\infty}$ have the
desired form.

Suppose that $m$ is a positive integer, $f=(f_1,\ldots,f_m)\in (V^*)^m$.
Set $u_i=y_i(0)$, $1\le i\le n$.
Set $u_{n+i}=f_i$ for $1\le i\le m$.
 If $f_i\in V^*\cap K=V$, also set $v_{i}=f_i$.

By Theorem 4.9 \cite{Dale} and Theorem 6.1 \cite{Dale}, with the
$R^*$, $S^*$ in the assumptions of Theorem 4.9 \cite{Dale} set as
$R^*=R_0(0)$, $S^*=S_0$, and with the $\{u_i\}$ and $\{v_i\}$
defined as above, and then applying Theorem \ref{TheoremFE6} (and
Remark \ref{RemarkFE19}), there exists a commutative diagram
$$
\begin{array}{lllll}
R_0(f)&\rightarrow&R_f&\rightarrow&S_f\\
\uparrow&&&&\uparrow\\
R_0(0)&\rightarrow&&&S_0
\end{array}
$$
such that the vertical arrows are sequences of monoidal transforms
along $V^*$, (1) and (2) of this theorem hold for
$$
R_0(f)\rightarrow R_f\rightarrow S_f
$$
and (3)(a) of this theorem holds for
$$
\begin{array}{lll}
R_f&\rightarrow&S_f\\
\uparrow&&\uparrow\\
R_0&\rightarrow&S_0
\end{array}.
$$
We have that $\text{det}(D(0,f))\ne0$ since (1) holds for $S_0$
and $S_f$. Define a partial ordering on $I=\sqcup_{m\in {\bold
N}_+}(V^*)^m$ by $f\le g$ if $S_f\subset S_g$. We will associate
to $0\in V^*$ the sequence $R_0(0)\rightarrow R_0\rightarrow S_0$
constructed in the beginning of the proof. Suppose that
$\alpha\le\beta$. We have
$$
S_0\subset S_{\alpha}\subset S_{\beta}
$$ so the proof of (2) of Theorem \ref{TheoremA3} shows that (3)(a) of this Theorem holds
for $\alpha,\beta$.

(3)(b) holds since
$$
R_{\alpha}=S_{\alpha}\cap K\subset S_{\beta}\cap K=R_{\beta}.
$$
(3)(c) is immediate, since the conclusions of Theorem 6.1
\cite{Dale} hold. In particular, (11) of Theorem 4.7 \cite{Dale}
holds.

Finally, we will establish (4) of the Theorem. By construction,
$V^*=\cup_{j\in I}S_j$. If $f\in V$, we have $f\in S_f\cap K=R_f$,
thus $V=\cup_{j\in I}R_j$. By construction, $\cup_{j\in
I}\Omega_j=\Gamma^*$, since $\nu^*(f)\in\Omega_f$ for $f\in V^*$.
We also have $\cup_{j\in I}\Lambda_j=\Gamma$, since
$\nu(f)\in\Lambda_f$ for $f\in V$. $I$ is a directed set as shown
in the proof of Theorem \ref{TheoremA3}.
\end{pf}

\begin{Theorem}\label{TheoremFE15} Let assumptions be as in Theorem \ref{Theorem8}.
There exists a partially ordered set $I$ and algebraic regular
local rings $\{S_i\mid i\in I\}$ with quotient field $K^*$ which
are dominated by $V^*$ and algebraic local rings with toric
singularities $\{R_i\mid i\in I\}$ such that
\begin{enumerate}
\item
$$
V^*=\lim_{\rightarrow} S_i=\cup_{i\in I}S_i, V=\lim_{\rightarrow}
R_i=\cup_{i\in I}R_i
$$
and each $S_i$ is a localization at a maximal ideal of the
integral closure of $R_i$ in $K^*$. \item $\hat S_j/p(\hat
S_j)_{\infty}$ are regular local rings for all $j$ and
$$
T^*=\lim_{\rightarrow}\hat S_j/p(\hat S_j)_{\infty}
$$
is a Henselian valuation ring such that $(T^*,Q(T^*))$  is an
immediate extension of $(V^*,K^*)$.

\item $\hat R_j/p(\hat R_j)_{\infty}$ has normal toric
singularities for all $j$ and
$$
T=\lim_{\rightarrow} \hat R_j/p(\hat R_j)_{\infty}
$$
 is a Henselian
valuation ring such that $(T,Q(T))$ is an immediate extension of $(V,K)$.

\item Further suppose that $k=V^*/m_{V^*}$ is algebraically closed
$($of characteristic zero$)$. Then the action of $\Gamma^*/\Gamma$
on $\hat S_j$ by $k$-algebra isomorphisms extends to an action of
$\Gamma^*/\Gamma$ on $\hat S_j/p(\hat S_j)_{\infty}$, and an
action of $\Gamma^*/\Gamma$ on $T^*$ such that $\hat R_j/p(\hat
R_j)_{\infty}\cong (\hat S_j/p(\hat
S_j)_{\infty})^{\Gamma^*/\Gamma}$ and $T=(T^*)^{\Gamma^*/\Gamma}$.
\end{enumerate}
\end{Theorem}
\begin{pf} The theorem follows from Theorems \ref{Theorem8}, \ref{TheoremA6}
 and \ref{TheoremA4}. The fact that $p(\hat S_j)_{\infty}$ is fixed by $\Gamma^*/\Gamma$
follows from (2) of Theorem \ref{Theorem8}.
\end{pf}

Now the proof of Theorem \ref{TheoremI3} follows from Theorem
\ref{TheoremFE15}. The proof that $Q(U^*)$ is Galois over $Q(U)$
is as in Theorem \ref{TheoremI2}.

\end{document}